\documentclass[a4paper,12pt,reqno]{amsart}

\usepackage[applemac]{inputenc}
\usepackage[T1]{fontenc}
\usepackage{graphicx}
\usepackage[toc,page]{appendix}
\usepackage{lscape}
\usepackage{amssymb, microtype,graphicx}
\usepackage{booktabs}  % For publication quality tables in LaTeX

\newcounter{notecounter}

\newtheorem{Th}{Theorem}

\newtheorem{Co}[Th]{Corollary}
\newtheorem{Ex}[Th]{Example}
\newtheorem{Le}[Th]{Lemma}

\theoremstyle{definition}

\title{The Farey maps modulo n}
\author{David Singerman and James Strudwick}

\date{\today}

\begin{document}
\maketitle
\footnote{University of Southampton}

\section{Abstract} The Farey map is the universal triangular map whose automorphism group is the classical modular group. We study the quotients of the Farey map by the principal congruence subgroups of the modular group. We also study the structure of the underlying graphs of these quotients.

\section{Introduction}
Let $\mathcal{U}^*$ denote the upper-half plane compactified by adding the points $\mathbb{Q}\cup\{\infty\}$ to the upper half plane $\mathcal{U}$.  On $\mathcal{U}^*$  we have the universal triangular map ${\mathcal {M}_3}$ which can be realised by the well-known Farey map as described below. The extended rationals, $\mathbb Q\cup\{\infty\}$, form the vertices of the map.
Our aim in this paper is to discuss the maps (or clean dessin d'enfants)  ${\mathcal{M}_3}/ \Gamma(n)$ which lie on the Riemann surface $\mathcal{U}^*/\Gamma(n)$ where $\Gamma(n)$ is the principal congruence subgroup mod  $n$ 
of the classical modular group $\Gamma$.  Their vertices can be identified with rational numbers "modulo $n$". Such maps were introduced in \cite{IS} and also discussed in \cite{SS}.

\section{The Farey map} 
 Vertices of the Farey map are the extended rationals, i.e. $\mathbb Q\cup\{\infty\}$ and two rationals $\frac{a}{c}$ and $\frac{b}{d}$ are joined by an edge if and only if $ad-bc=\pm 1$. These edges are drawn as semicircles or vertical lines, perpendicular to the real axis, (i.e. hyperbolic lines).
 Here $\infty=\frac{1}{0}$. This map has the following properties. 
\begin{enumerate}
\item[(a)] There is a triangle with vertices  $\frac{1}{0}, \frac{1}{1}, \frac{0}{1}$, called the principal triangle.
\item[(b)] The modular group $\Gamma$= PSL(2, $\mathbb{Z})$ acts as a group of automorphisms of $\mathcal{M}_3$.
\item[(c)] The general triangle has vertices $\frac{a}{c},\frac{a+b}{c+d}, \frac{b}{d}$.
\end{enumerate}
This forms a triangular tessellation of the upper half plane. Note that the triangle in (c) is just the image of the principal triangle under the M\"obius transformation corresponding to the matrix $\begin{pmatrix} 
  a&b\\
  c&d\\
  \end{pmatrix}$.

In \cite {S} it was shown that ${\mathcal {M}_3}$  is the universal triangular map in the sense that any triangular map on a surface is the quotient of ${\mathcal {M}_3}$ by a subgroup of the modular group $\Gamma$ and any regular triangular map is the quotient of ${\mathcal {M}_3}$ by a normal subgroup of $\Gamma$. The aim of this paper is to study the maps $\mathcal{M}_3(n)=\mathcal{M}_3/\Gamma(n)$ where $\Gamma(n)$ is the principal congruence subgroup mod $n$ of the Modular group $\Gamma$, see \S5.

\begin{figure}
\begin{center}
\includegraphics[height=\linewidth,angle=-90]{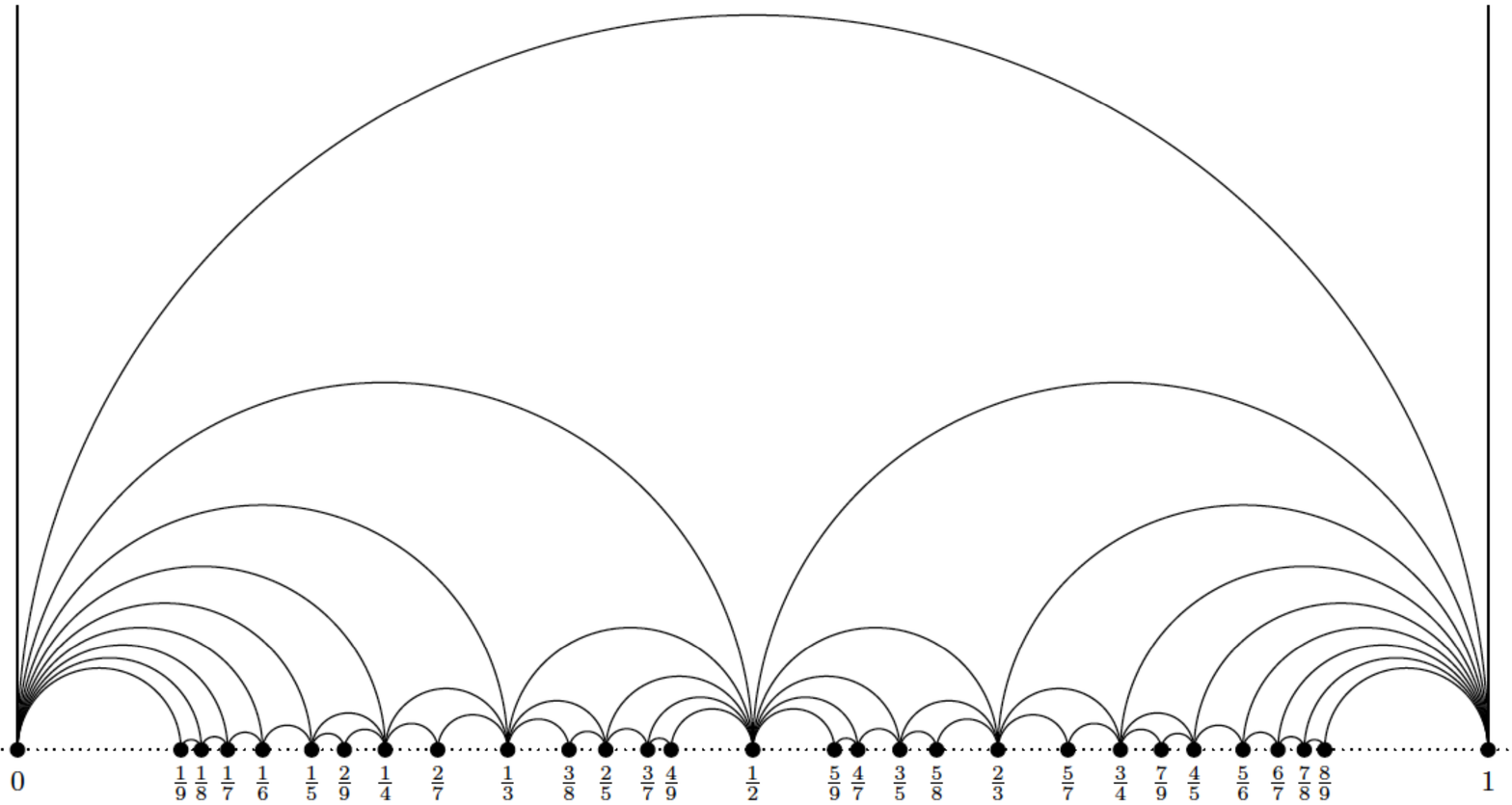}
\end{center}
\caption{The Farey map, (drawn by Jan Karaba\v s)}
\end{figure}
 
\section{The graph underlying $\mathcal{M}_3(n)$} 
The vertices of this graph are the Farey rationals mod $n$.  These rationals are of the form $\frac{a}{c}$ where $(a,c,n)=1$, excluding $\frac{0}{0}$. We can think of these as fractions $\frac{a}{c}$ where now $a,c\in \mathbb{Z}_n$, (not both 0) and where we identify $\frac{a}{c}$ with $\frac{-a}{-c}$ and two vertices $\frac{a}{b}$ and $\frac{c}{d}$ are joined by an edge if and only if $ad-bc\equiv\pm 1 \pmod n$.
 
\subsection{Low values of $n$}   
Let us see what these vertices are for low values of $n$, say $n=1,2\cdots, 8.$ As we will see, these maps give interesting geometric objects. See Appendix A.

$n=2$. The vertices are $\frac{1}{0},\frac{0}{1}, \frac{1}{1}$.

We are here abusing notation, by $\frac{0}{1}$ we are regarding $0$ and $1$ as congruence classes mod 2 as well as integers. This abuse of notation will be used throughout. It should not cause confusion.

$n=3$. The vertices are $\frac{1}{0},\frac{0}{1},\frac{1}{1},\frac{2}{1}.$

$n=4$. The vertices are $\frac{1}{0}, \frac{0}{1}, \frac{1}{1}, \frac{2}{1},\frac{3}{1}, \frac{1}{2}.$

$n=5$. The vertices are $\frac{1}{0}, \frac{2}{0}, \frac{0}{1}, \frac{1}{1}, \frac{2}{1}, \frac{3}{1},\frac{4}{1},\frac{0}{2},\frac{1}{2}, \frac{2}{2},\frac{3}{2}, \frac{4}{2}.$

$n=6$. The vertices are  $\frac{1}{0}, \frac{0}{1}, \frac{1}{1},  \frac{2}{1},\frac{3}{1}, \frac{4}{1}, \frac{5}{1}, \frac{1}{2}, \frac{3}{2}, \frac{5}{2},\frac{1}{3},\frac{2}{3}.$

$n=7$. The vertices are  $\frac{1}{0}, \frac{2}{0}, \frac{3}{0}, \frac{0}{1}, \frac{1}{1},  \frac{2}{1},\frac{3}{1}, \frac{4}{1}, \frac{5}{1},\frac{6}{1}, \frac{0}{2}, \frac{1}{2}, \frac{2}{2},\frac{3}{2}, \frac{4}{2},\frac{5}{2},\frac{6}{2},\frac{0}{3}, \frac{1}{3}, \frac{2}{3},  \frac{3}{3},\frac{4}{3},\frac{5}{3},\frac{6}{3}.$

$n=8$. The vertices are $\frac{1}{0}, \frac{3}{0},\frac{0}{1},\frac{1}{1}, \frac{2}{1},\frac{3}{1}, \frac{4}{1}, \frac{5}{1}, \frac{6}{1},\frac{7}{1}, \frac{1}{2},\frac{3}{2}, \frac{5}{2}, \frac{7}{2},\frac{0}{3}, \frac{1}{3}, \frac{2}{3}.\frac{3}{3}, \frac{4}{3}, \frac{5}{3},\frac{6}{3},\frac{7}{3}.
\frac{1}{4},\frac{3}{4}.$

 \section{PSL(2, $\mathbb{Z}_n)$ as a $(2,m,n)$ group} In this section we go over some well-known material which some may wish to ignore. 
 A  group $G$ is called a $(2,m,n)$-group if is a homomorphic image of the $(2,m,n)$ triangle group $\Gamma(2,m,n)$. This is the group with generators $(X,Y,Z)$ which obey the relations $X^2=Y^m=Z^n=XYZ=1$. If $G$ is a finite $(2,m,n)$-group then there is an epimorphism $\theta:\Gamma(2,m,n)\mapsto G$ with kernel $M$, say.  If $\mathcal{U}$ is the upper-half complex plane then $\mathcal{U}/M$ is a compact Riemann surface that carries a map $\mathcal{M}$ of type $(m,n)$.  In \cite{JS}, $M$ is called a {\it map subgroup} for $\mathcal{M}$. If we are just interested in triangular maps then we only deal with $(2,3,n)$ triangle groups. If we are not particularly concerned with the vertex valencies then we are just looking at factor groups of the $(2,3,\infty)$ triangle group which means that we are not concerned with the order of the image of $Z$.  Now the $(2,3,\infty)$ triangle group is isomorphic to the classical modular group  $\Gamma$, which is defined by: 
    $$\Gamma=\Big\{\begin{pmatrix} 
  a&b\\
  c&d\\
  \end{pmatrix}: a,b,c,d\in\mathbb{Z}, \quad ad-bc=\pm1\Big\}/\{\pm I\}$$. 
 
  The definition of the principal congruence subgroup $\Gamma(n)$ is as follows   
    $$\Gamma(n)=\Big\{\begin{pmatrix} 
  a&b\\
  c&d\\
  \end{pmatrix}\in \Gamma: 
  \begin{pmatrix}
  a&b\\
  c&d\\
\end{pmatrix}\equiv
\pm\begin{pmatrix}1&0\\0&1
\end{pmatrix}\text {mod}\;n\Big\}$$ 
  
 %It is well-known that $\Gamma$ is generated by 
%\begin{center} 
 %$X=\begin{pmatrix}0&1\\-1&0\end{pmatrix}$, $Y=\begin{pmatrix}0&1\\-1&1\end{pmatrix}$
 %\end{center}
 %(As usual these matrices represent M\"obius transformations, so we are working in $PSL(2,\mathbb{Z})=SL(2,\mathbb{Z})/\{\pm 1\}$.) Let 
 %\begin{center}
%$Z=XY=\begin{pmatrix}-1&1\\0&-1\end{pmatrix}$  
 %\end{center}
%which has infinite order in $\Gamma$. 
Now $\Gamma(n)\triangleleft \Gamma$ and the quotient group $\Gamma/\Gamma(n)=PSL(2,\mathbb{Z}_n)$.  This group has the same generators as $\Gamma$ except now the integers are taken modulo $n$ and hence $X^2=Y^3=Z^n=XYZ=1$ and thus $PSL(2,\mathbb{Z}_n)$ is a finite $(2,3,n)$ group.   There is an epimorphism $\theta:\Gamma\rightarrow PSL(2, \mathbb{Z}_n)$. The kernel of $\theta$ is denoted by $\Gamma(n)$ and called the prinicipal congruence subgroup of $\Gamma$ of level $n$.  Now $\Gamma(n)$ is a map subgroup for a regular triangular map which we denote by $\mathcal{M}_3(n)$.

\section{The vital statistics of  ${\mathcal {M}_3}(n)$}

  Every element of $\Gamma$ is represented by a matrix $\begin{pmatrix} a&b\\ c&d\\\end{pmatrix}$. Now $\Gamma$ acts transitively on  the darts of $\mathcal{M}_3$, where a {\it dart} is a directed edge. A dart of $\mathcal{M}_3$ is an ordered pair ($\frac{a}{c}, \frac{b}{d}$) and the matrix represented by $\begin{pmatrix}a&b\\ c&d\\\end{pmatrix}$  maps the dart ($\frac{1}{0}, \frac{0}{1})$ to  $(\frac{a}{c}, \frac{b}{d})$ and so $G=\Gamma/\Gamma(n)$ acts transitvely on  the darts of $\mathcal{M}_3/\Gamma(n)$. By definition, a map is regular if its automorphism group acts transitively on its darts which makes  $\mathcal{M}_3$ a regular map and for similar reasons so is $\mathcal{M}_3(n)$. Note that $\frac{1}{0}$ is joined to $\frac{k}{1}$   in $\mathcal{M}_3(n)$ for $k=0,1,2,\cdots, \frac{n-1}{1}$, so that $\frac{1}{0}$ has valency $n$ and hence by regularity, every vertex has valency $n$. We want to find the numbers of darts, edges, faces and vertices of $\mathcal {M}_3(n)$.  Being a regular map, for $\mathcal{M}_3(n)$, we have:
  \begin{enumerate}
  \item[-] The number of darts = $|G|$
  \item[-] The number of edges = $|G|/2$
  \item[-] The number of faces = $|G|/3$ 
  \item[-] Every vertex has valency $n$
  \item[-] The number of vertices = $|G|/n$
  \end{enumerate}
%The definition of the principal congruence subgroup $\Gamma(n)$ is as follows  
 % $$\Gamma(n)=\Big\{\begin{pmatrix} 
  %a&b\\
  %c&d\\
  %\end{pmatrix}\in \Gamma: 
  %\begin{pmatrix}
  %a&b\\
  %c&d\\
%\end{pmatrix}\equiv
%\pm\begin{pmatrix}1&0\\0&1
%\end{pmatrix}\text {mod}\;n\Big\}$$
  
    For $n>2, \Gamma(n)$ is a normal subgroup of $\Gamma$ of index
 $$\mu(n)=\frac{n^3}{2} \Pi_{p|n}(1-\frac{1}{p^2}). \eqno(1)$$ For $n=2$ the index is equal to 6.
 For this see any text on modular functions, e.g. \cite{DS}.
 Thus $|G|=\mu(n)$ is the number of darts so that the number of edges, vertices and faces is $\mu(n)/2$, $\mu(n)/n$ and $\mu(n)/3.$
 If $g(n)$ is the genus of the map $\mathcal{M}_3(n)$ then the Euler characteristic is given by 
 
 $$2-2g(n)=\mu(n)(\frac{1}{n}-\frac{1}{2}+\frac{1}{3})=\mu(n)\frac{(6-n)}{6n}$$
 
 from which we deduce the following formula for the genus of $\mathcal{M}_3(n)$
 
 $$g(n)=1+\frac{n^2}{24}(n-6)\Pi_{p|n}(1-\frac{1}{p^2})$$
 
 which is the well-known formula for the genus of $\Gamma(n)$, i.e. the genus of the surface $\mathcal{U}/\Gamma(n)$ which carries the map $\mathcal{M}_3(n)$. 

We can also find the number of vertices using group theory.  An important subgroup of $\Gamma$ is 
 
 $$\Gamma_1(n)=\Big\{\begin{pmatrix} 
  a&b\\
  c&d\\
  \end{pmatrix}\in \Gamma: 
  \begin{pmatrix}
  a&b\\
  c&d\\
\end{pmatrix}\equiv 
\pm\begin{pmatrix}1&b\\0&1
\end{pmatrix}\text {mod}\; n\Big\}$$

where $0\le b<n$. As the stabilizer of $\infty$ in $\Gamma$ is just the infinite cyclic group generated by $U=\begin{pmatrix}1&1\\0&1\end{pmatrix}$
it easily follows that the stabilizer of $[\infty]_\Gamma(n)\in {\mathcal {M}_3}(n)$ is equal to $\Gamma_1(n)/\Gamma(n)$.
(See Lemma 3.2 in \cite{IS}). By the Orbit-Stabilizer Theorem  there exists a one-to one correspondence between the left cosets of $\Gamma_1(n)$ in $\Gamma$ and the 
vertices of $\mathcal{M}_3(n)$.

Now there is a epimomorphism $\chi:\Gamma_1(n)\mapsto \mathbb{Z}_n$ defined by 
$$\chi\Big(\pm\begin{pmatrix} an+1&b\\cn&dn+1\end{pmatrix}\Big)\mapsto b(mod {\text\;n}).$$
The kernel of $\chi$ is $\Gamma(n)$. For $n>2, \,  \Gamma_1(n)$ is a normal subgroup of $\Gamma$ of index 
$$\frac{n^2}{2} \Pi_{p|n}(1-\frac{1}{p^2}) \eqno(2)$$ and 
the index of $\Gamma(2)$ in $\Gamma$ is equal to 6, Hence  the number of vertices of ${\mathcal {M}_3}(n)$ is given by (2) if $n>2$ and the number of vertices of ${\mathcal{M}_3}(2)$ is equal to 3. The reader is invited to check this formula for $n=2\dots8$ by the lists given in Section 3.1.

\begin{Ex}  The maps ${\mathcal {M}_3}(n)$ were described in \cite{IS} for $n=2,...,7$. For $n=2$ we get a triangle, for $n=3$ we get a tetrahedron, for $n=4$ we get an octahedron and for $n=5$ we get an icosahedron. All these are planar maps and so the corresponding Riemann surface is the Riemann sphere. For $n=6$ we get the torus map $\{3,6\}_{2,2}$ and the corresponding Riemann surface is the hexagonal torus and for $n=7$ we get the Klein map which lies on Klein's Riemann surface of genus 3, known as the Klein quartic \cite{K}. See figures 2,3,4 in the appendix.
\end{Ex}
 
\begin{Ex}{$n=8$ } We find that $\mu(8)=192$, so that the number of edges of $\mathcal{M}_3(8)$ is equal to 96, the number of faces is 64 and the number of vertices is 24 and so the genus of $\mathcal{M}_3(8)$ is equal to 5. The underlying Riemann surface is the unique Riemann surface of genus 5 with 192 automorphisms and this is known as the {\it Fricke-Klein surface of genus 5}. This map is particularly important as it is the regular map underlying the \emph {Grunbaum polyhedron} which is one of the few regular  maps that are known to admit a polyhedral embedding into Euclidean 3-space with convex faces \cite{GSW}.  For the uniqueness of this map and Riemann surface we refer to \cite{CD}. This map and Riemann surface also appears in \cite{I}. See figure 5 in the appendix.
\end{Ex}
%For pictures of all these maps see the appendix at the end of this paper.

\section{Petrie paths in $\mathcal{M}_3(n)$}
A {\it Petrie path} in a (regular) map is a zig-zag path on the map. This is a path in which two consecutive edges but no three consecutive edges can belong to the same face. If this is a closed path then it is called a \emph{Petrie polygon} and the number of edges of this polygon is called the {\it Petrie length} of the map.

In \cite{SS} it was shown that the $k$th vertex of the Petrie path of the universal map is $\frac{f_{k-1}}{f_k}$ where $f_k$ is the $k$th element of the Fibonacci sequence given by $f_0=1, f_1=0, f_{k+1}=f_k+f_{k-1}.$. For $\mathcal{M}_3(n)$ we get exactly the same results except now the integers $f_i$ are taken modulo $n$.  As an example we give the result from \cite{SS}.  The Petrie polygon in $\mathcal{M}_3(7)$ is $\frac{1}{0},\frac{0}{1},\frac{1}{1},\frac{1}{2}, \frac{2}{3}, \frac{3}{5}, \frac{5}{1},\frac{1}{6}, \frac{6}{0}=\frac{-1}{0}$ so we have closed the path which in this case has length equal to 8. In general, we find that the length of the Petrie polygon on $\mathcal{M}_3(n)$ is equal to $\sigma(n)$ where $\sigma(n)$ is the {\it semiperiod} of the Fibonacci sequence mod $n$. This means the period of the Fibonacci sequence up to sign, so for $n=7$ we get 1,0,1,1,2,3,5,1,-1, 0,-1,-1..so  $\sigma(7)=8$. (This was one of the main topics of \cite{SS}.) As another example, let us find the Petrie path in $\mathcal{M}_3(8).$   The reader will enjoy finding this and it is seen that it has length $\sigma(8)=12.$ See figure 5 in the appendix noting that the Petrie paths are drawn in red.

 \section{The star of a vertex}
  In a graph $G$ the star of a vertex $x$ consists of all vertices of $G$ that are joined to $x$ by an edge, including $x$ itself. Now let $\frac{a}{c}$ be a vertex of ${\mathcal {M}_3}(n)$, where $(a, c, n)=1$. 
  
\begin{Le}  There exists $b,d\in \mathbb{Z}$ such that $ad-bc\equiv 1\pmod n$.  
\begin{proof} Suppose that $(a,c)=k$.  Then $(k,n)=1$ so there exist $\alpha,\beta\in \mathbb{Z}$ such that $\alpha k+\beta n=1$. Now there exists $u,v\in \mathbb{Z}$ such that $ua+vc=k$ so that $\alpha(ua+vc)+\beta n=1$ and hence $(\alpha u) a+(\alpha v)c +\beta n=1$.
\end{proof}
\end{Le}
  
\begin{Th}  The star of $\frac{a}{c}$ consists of $\frac{a}{c}$ together with all vertices of the form $\frac{ak+b}{ck+d}$ where $k=0,1,\cdots ,n-1$ and $b,d$ are as in lemma 3.
\begin{proof}  First we find the star of $\frac{1}{0}$. This consists of $\frac{1}{0}$ together with the vertices$\frac{0}{1}, \frac{1}{1}, \cdots \frac{n-1}{1}$. Here $\frac{a}{c}=\frac{1}{0}$ so $\frac{b}{d}=\frac{0}{1}$ so $a=d=1$ and $b=c=0$ so $\frac{b+ak}{d+ck}=\frac{k}{1}$, $k=0,1,\dots, {n-1}$ as required. More generally, let
$T=\begin{pmatrix}a&b\\c&d\end{pmatrix}$, $U=\begin{pmatrix}1&1\\0&1\end{pmatrix}$.  Then $T(\frac{1}{0})=\frac{a}{c}$ . The stabilizer of $\frac{1}{0}$ is the cyclic group generated by $U$ so the stabilizer of $\frac{a}{c}$ consists of elements of the form  by $TU^kT^{-1}$.  Now 
$$S:=TU^kT^{-1}=\begin{pmatrix}1-akc&a^2k\\-c^2k&1+akc\end{pmatrix}$$
Therefore if $\frac{a}{c} \longleftrightarrow \frac{b}{d}$ then $S(\frac{a}{c}) \longleftrightarrow S(\frac{b}{d})$ so that $\frac{a}{c}\longleftrightarrow S(\frac{b}{d})$.
Now $$S\left(\frac{b}{d}\right)=\begin{pmatrix}1-akc&a^2\\-c^2k&1+akc\end{pmatrix}\left(\frac{b}{d}\right)$$
$$=\frac{b-abck+a^2kd}{-bc^2k+d+acdk}$$
On using $ad=1+bc$ this becomes
$$\frac{ak+b}{ck+d}.  \eqno(3)$$
This is true for $k=0,1,\cdots, (n-1)$ so that these are the $n$ points in the star of $\frac{a}{c}$.
\end{proof}
\end{Th}
  
Thus the numerators form an arithmetic progression of length $n$ whose first term is $b$ and common difference is $a$ and the denominators form an arithmetic progression of length $n$ whose first term is $d$ and common difference is $c$.
 
Therefore to find the star of $\frac{a}{c}$ we find $b,d\in \mathbb{Z}$ such that $ad-bc\equiv 1 \pmod n$, construct the unimodular matrix $\begin{pmatrix}a&b\\c&d\end{pmatrix}$ and then calculate the star as is given by (3).
 
\begin{Ex} Find the star of $\frac{3}{5}$ in $\mathcal {M}_3(7)$.  Our unimodular matrix is now $\begin{pmatrix}3&0\\5&5\end{pmatrix}$ so that the numerators form an arithmetic progression whose first term is 0 and common difference 3 and whose denominators form an arithmetic progression whose first term is 5 and whose common difference is 5.  Thus the star of $\frac{3}{5}$ in $\mathcal {M}_3(7)$ is $\frac{3}{5}$ together with
$$\left\{\frac{0}{5}, \frac{3}{3}, \frac{6}{1}, \frac{2}{6}, \frac{5}{4}, \frac{1}{2}, \frac{4}{0}\right\}.$$ 
\end{Ex}
 
See Figure 4 noting that $\frac{3}{5}=\frac{4}{2}$.

 \begin{Ex}Find the star of $\frac{1}{2}$ in $\mathcal {M}_3(8)$.  Our unimodular matrix is now $\begin{pmatrix}1&2\\2&3\end{pmatrix}$ so that the numerators form an arithmetic progression whose first term is 2 and whose common difference is 1 and the denominators form an arithmetic progression whose first term is 3 and whose common difference is 2.
 Thus the star of $\frac{1}{2}$ in $\mathcal {M}_3(8)$ is $\frac{1}{2}$ together with
$$\left\{\frac{2}{3},  \frac{3}{5},  \frac{4}{7},  \frac{5}{1}, \frac{6}{3}, \frac{7}{5}, \frac{0}{7} ,\frac{1}{1}\right\}.$$ 

 \end{Ex}
 
See figure 5 in the appendix, noting that $\frac{3}{5}=\frac{5}{3}$, $\frac{4}{7}=\frac{4}{1}$ and $\frac{7}{5}=\frac{1}{3}$. Note that these give the cyclic orderings of the stars around our two points $\frac{3}{5}$ and $\frac{1}{2}$.  This is because the star of $\frac{1}{0}$ is $\frac{0}{1},\frac{1}{1},\frac{2}{1}\dots$.
 
 {\bf Note} The vertices of the star of $\frac{a}{c}$ form a polygonal face whose centre is $\frac{a}{c}$.
 \newpage
 \section{The stars of $\mathcal{M}_3(p)$, $p$ an odd prime}
 
\begin{Th} For a prime, $p$, the stars of $\frac{1}{0}, \frac{2}{0}, \dots , \frac{(p-1)/2}{0}$ are disjoint and cover $\mathcal{M}_3(p)$.
\begin{proof} 
Consider the vertex $\frac{k}{0}$.  Let $K$ be the inverse of $k$ mod $p$. Then the star of $\frac{k}{0}$ is $\frac{k}{0}$ together with
$$\left\{ \frac{0}{K},\frac{1}{K}\dots, \frac{(p-1)}{K}\right\}.$$
Consider a distinct vertex $\frac{l}{0}$, so that $l\not= \pm k$.  Let $L$ be the inverse of $l$ mod $p$.  Then $L \not=\pm K$ so that the star of $\frac{l}{0}$ is $\frac{l}{0}$ together with
$$\left\{\frac{0}{L}, \frac{1}{L},\dots,   \frac{(p-1)}{L}\right\}.$$
As $L\not=\pm K$ these stars are disjoint.  There are $(p-1)/2)$ stars each containing $p$ vertices. Thus there are $p(p-1)/2$ vertices which is the total number of vertices $\mathcal{M}_3(p)$. \end{proof}
\end{Th}

\begin{Ex} $\mathcal{M}_3(5).$  We find the stars of $\frac{1}{0}$ and $\frac{2}{0}$.
The star of $\frac{1}{0}$ $\frac{1}{0}$ together with is $\{\frac{0}{1}, \frac{1}{1}, \frac{2}{1}, \frac{3}{1}, \frac{4}{1}.\}$.  
The star of $\frac{2}{0}$ is $\frac{2}{0}$ together with $\{\frac{0}{3},\frac{2}{3}, \frac{4}{3}, \frac{1}{3}, \frac{3}{3}\}$.  Thus the stars of $\frac{1}{0}$ and $\frac{2}{0}$ give us 12 vertices which is the number of vertices of  $\mathcal{M}_3(5)$,

\end{Ex}
  
 % \begin{figure}
%\begin{center}
%\includegraphics[width=\linewidth]{M5}
%$\end{center}
%\caption{Icosahedron}
%\end{figure}
\section{Poles of Farey maps}
  From the diagrams in the appendix we see that the points of the form $\frac{a}{0}$ play a significant role. We call these the \emph {poles} of $\mathcal{M}_3(n)$. As $\frac{x}{y}=\frac{-x}{-y}$ we see that for $n>2$ the number of poles of $\mathcal{M}_3(n)$ is equal to $\phi(n)/2$.  %This is equal to 1 if and only if $\phi(n)=2$ corresponding to the line, triangle, octahedron, or the toroidal map $\{3,6\}_{2,2}$ in the appendix. Perhaps more interesting is the case when there are two poles. Now $\phi(n)=4$ so that $n=5,8,10, 12.$ The first two cases correspond to $n=5$ and $n=8$, the icosahedron and the Fricke-Klein map. These are mentioned in \S 4 and drawn in Figures. 2 and 5 in the appendix.

\subsection{Graphical distance} If $u$,$v$ are two vertices in a graph then the distance $\delta(u,v)$ between them is defined to be the length of the shortest path joining $u$ and $v$,  The {\it diameter} of a graph or map is the maximum distance between two of its points.

\begin{Th}\label{dis3} 
The distance between two distinct poles of $\mathcal{M}_3(n)$ is equal to 3.
\begin{proof} By regularity we may assume that one of the poles is $\frac{1}{0}$. This is not adjacent to $\frac{a}{0}$.  There is no path of length 2 between $\frac{1}{0}$ and $\frac{a}{0}$, where $a\not=\pm 1$, for otherwise there would be $x, y\in \mathbb{Z}$ such that $\frac{1}{0}\rightarrow \frac{x}{y}\rightarrow\frac{a}{0}$.  Then $y=\pm 1$ and then $a=\pm 1$, a contradiction. However, we can always construct a path of length 3 of the form 
$$\frac{1}{0}\rightarrow \frac{x}{1}\rightarrow \frac{xa^{-1}+1}{a^{-1}}\rightarrow \frac{a}{0}.$$
where $a^{-1}$ is the inverse of $a$ modulo $n$. This inverse exists because $(a,0,n)=1$ and hence $(a,n)=1$.
\end{proof}
\end{Th}

\begin{Le} Let $a,c,n$ be integers such that $(a,c,n)=1$. Then there exists an integer $k$ so that $(a+ck, n)=1$. \label{unitProof}

\begin{proof} Let $p_1, p_2\dots ,p_r$ be a list without repetition of those prime divisors of $n$ that are not divisors of $a$ or $c$. . Define $k=p_1\cdot.....p_r$, the product of the $p_i$.

Now suppose that $q$ is a prime divisor of $n$. If $q=p_i$ for some $i$ then $q$ is not a divisor of $a$ so is not a divisor of $a+ck$.

If $q$ is not equal to $p_i$ for any $i$, then it is coprime to $k$ and is a divisor of one, but not both of $a$ and $c$. If $q$ were a divisor of both $a$ and $c$ then $(a,c,n)=q\not=1$.  If $q$ is a divisor of $a$, then it cannot divide $a+ck$, as $q$ is coprime to $ck$.  If $q$ is a divisor of $c$, then $q$ does not divide $a$ and so does not divide $a+ck$.
Hence $a+ck$ and $n$ are coprime..

\end{proof}
\end{Le}

\begin{Th} The diameter of $\mathcal{M}_3(n)$ is equal to $3$ for all $n\geq5$. \label{diam}
\begin{proof}
By transitivity we can assume one of the points to be $\frac{1}{0}$ and denote the other by $\frac{a}{c}$. Now the neighbours of $\frac{1}{0}$ are precisely those of the form $\frac{k}{1}$ for $k=0,\dots,n-1$. This reduces to finding an $\frac{x}{y}\in\mathcal{M}_3(n)$ such that $\frac{k}{1}\rightarrow\frac{x}{y}\rightarrow\frac{a}{c}$.

In terms of congruences this means that we have two simultaneous equations:
\begin{align*}
ky-x & \equiv 1 \mod(n)\\
cx-ay & \equiv 1 \mod(n)\\
\Rightarrow y(ck-a) & \equiv c+1 \mod(n)\label{finEqu}
\end{align*}
By Lemma \ref{unitProof} we know that $(ck-a)$ is co-prime to $n$ and therefore has a inverse mod $n$. Hence the last equation above can be solved for $y$ which, in turn, determines $x$.
\end{proof}
\end{Th}

For $n=2,3,4$ the diameters are $1,1,2$ respectively. For $n=2$ this corresponds to the triangle and for $n=3,4$ this is the tetrahedron and octahedron respectively which can be observed in Figure \ref{M345}.

\begin{Th} \label{dis2}
Given a point in $\mathcal{M}_3(n)$, $\frac{b}{d}$ where $d\not=\pm1$ then $\delta(\frac{1}{0},\frac{b}{d})=2$ if and only if $gcd(d,n)|(b\pm1)$.
\begin{proof}
By our supposition $\frac{1}{0}$ is not adjacent to $\frac{b}{d}$ therefore assume that there exist a vertex $\frac{x}{y}$ such that $\frac{1}{0}\rightarrow\frac{x}{y}\rightarrow\frac{b}{d}$. Now for $\frac{x}{y}$ to be adjacent to $\frac{1}{0}$ we must have $y=\pm1$. Setting $y=1$ and looking at $\frac{x}{y}\rightarrow\frac{b}{d}$ we arrive at the equation:
\begin{displaymath}
dx-b \equiv \pm1\mod(n)
\end{displaymath} 
which has solutions if and only if $gcd(d,n)|(b\pm1)$.
\end{proof}
\end{Th}

To see if any two points, $\frac{a}{c}$, $\frac{x}{y}$, are distance two apart let $T$ be the transformation where $T(\frac{a}{c})=\frac{1}{0}$. Apply this transformation to the other vertex, $T(\frac{x}{y})=\frac{b}{d}$, and then apply Theorem \ref{dis2}.

\begin{Co}\label{dis2prime}
Given a point $\frac{b}{d}$ in $\mathcal{M}_3(p)$, where $p$ is a prime. Then $\delta(\frac{1}{0},\frac{b}{d})=2$ if and only if $d\not=0,\pm1$. 
\begin{proof}
From Theorem \ref{dis2} to show $\delta(\frac{1}{0},\frac{b}{d})=2$, where $d\not=\pm1$, we had to solve equation:
\begin{displaymath}
dx \equiv b \pm1\mod(p)
\end{displaymath}
As we are now working mod$(p)$ this can be solved $\forall d\not=0$ as $gcd(d,p)=1$. If $d=0$ this would imply that $b=\pm1$, resulting in both points being equal.
\end{proof}
\end{Co}

\begin{Th}\label{dis3prime}
Given two points, $\frac{a}{c}, \frac{b}{d}$ in $\mathcal{M}_3(p)$, where $p$ is prime, then  $\delta\left(\frac{a}{c},\frac{b}{d}\right)=3$ if and only if $\Delta=0$. Where $\Delta=ad-bc$.
\begin{proof}
\begin{enumerate}
\item[($\Rightarrow$)]Apply the transformation $T\in\Gamma$ such that $T(\frac{a}{c})=\frac{1}{0}$ and $T(\frac{b}{d})=\frac{e}{f}$.  Now consider these two points. We have that $\delta(\frac{1}{0},\frac{e}{f})=3$ if and only if $f=0$. By definition they are adjacent if and only if $f=\pm1$. By Corollary \ref{dis2prime} the distance would be $2$ if and only if $f\not=0,\pm1$. In terms of residues this leaves only the case when $f=0$. From Theorem \ref{diam} the maximum distance is $3$ which therefore can only be the case when $f=0$. The determinant of these two points, $\frac{1}{0}$ and $\frac{e}{0}$, is equal to $0$. As $T$ preserves the determinants of points we see that the determinants of the original points is also equal zero.\\

\item[($\Leftarrow$)]  If $\Delta=0$ then we can find a $T\in \Gamma$ such that $T\left(\frac{a}{c},\frac{b}{d}\right)= \left(\frac{x}{0},\frac{y}{0}\right).$ 
Let $r=c, s=-a$. Then $ra+sc=0$. Also $(r,s,p)=1$ so there exist $u,q\in \mathbb{Z}$ such that $us-qr=1\mod(p)$ so let 
$T=\begin{pmatrix} 
  u&q\\
 r&s\\
\end{pmatrix}$
Then $T(\frac{a}{c})=\frac{ua+qc}{0}, T(\frac{b}{d})=\frac{ub+qd}{0}$. Hence $\frac{a}{c}$ and $\frac{b}{d}$ are in the orbit of poles and therefore, from Theorem \ref{dis3}, have $\delta(\frac{a}{c},\frac{b}{d})=3$. 
\end{enumerate}
\end{proof}
\end{Th}

Therefore we can now completely categorize the distances in $\mathcal{M}_3(p)$ when $p$ is prime 
\newpage

\begin{Th} Let $\frac{a}{c}, \frac{b}{d}$ be distinct Farey fractions in $\mathcal{M}_3(p)$, where $p$ is prime, and let $\Delta=ad-bc$.  Then:  
\begin{displaymath}
\delta\left(\frac{a}{c},\frac{b}{d}\right)=
\begin{cases}
1 & \mbox{if and only if $|\Delta|=1$,}\\
3 & \mbox{if and only if $\Delta=0$,}\\
2 & \mbox{otherwise.}
\end{cases} 
\end{displaymath}
\begin{proof} ($1$) Follows from definition and (2) follows from Theorem \ref{dis3prime}. For part (3) we know from Theorem \ref{diam} that the maximum distance between two distinct points is $3$. But in the first two parts when have shown the conditions for distance $1$ and $3$ are $\Delta=\pm1$ and $\Delta=0$. The only remaining option for points with $\Delta\not=0,\pm1$ is to be distance $2$ apart.
\end{proof}
\end{Th}

\section{Summary}
%To surmise in this paper we have continued the study of principal congruence maps, $\mathcal{M}_3(n)$, focusing on the arithmetic structure that arises within them. We have shown that the neighbors of a given vertex form an arithmetic progression determined by the vertex. In particular we showed that for when $n$ is prime these maps can in fact be decomposed into disjoint stars. We have also proved that the surfaces obey a sense of compactness in that for all $n\geq5$ the maximum distance between any two points is $3$. This maximum distance is achieved by, but not limited to, any unique pair of poles. Finally for prime $n$ we are able to determine the distance between any two points from the determinant they omit. 

We have shown that $\mathcal{M}_3(n)$, the Farey maps modulo $n$, even though they have a simple arithmetic definition, give rise to many interesting maps and Riemann Surfaces. In Section 9 we show that simple graph-theoretic properties have a simple connection to arithmetic progressions modulo $n$, (Theorem 4). We showed that for when $n$ is prime these maps can in fact be decomposed into disjoint stars, (Theorem 7). In Section 10 we discussed the poles of the Farey map; these include the centre point $\frac{1}{0}$ but otherwise they are points furthest away from the centre (Theorem 9). In Section 10 we also discussed graphical distance and showed in Theorem 11, that even though these maps can be quite big ($\mathcal{M}_3(n)$ has about $n^2$ verticies) they all have diameter $3$ for $n\geq5$. Finally for prime $n$ we are able to arithmetically charactise the distance between any two points (Theorem 15).

\bigskip

We would like to thank Ian Short for Lemma \ref{unitProof}, Gareth Jones, Juergen Wolfart and the referees for carefully reading this manuscript and suggesting some improvements.

\bigskip

.

  \newpage
  
 \bibliography{refFile}

\begin{thebibliography}{1}

\bibitem{CD}
{\sc Conder, M., and Dobcs{\'a}nyi, P.}
\newblock Determination of all regular maps of small genus.
\newblock {\em Journal of Combinatorial Theory, Series B 81}, 2 (2001),
  224--242.

\bibitem{DS}
{\sc Diamond, F., and Shurman, J.~M.}
\newblock {\em A first course in modular forms}, vol.~228.
\newblock Springer, 2005.

\bibitem{GSW}
{\sc G{\'e}vay, G., Schulte, E., and Wills, J.~M.}
\newblock The regular {G}r{\"u}nbaum polyhedron of genus 5.
\newblock {\em Advances in Geometry 14}, 3 (2014), 465--482.

\bibitem{I}
{\sc Ivrissimtzis, I., Peyerimhoff, N., and Vdovina, A.}
\newblock Trivalent expanders and hyperbolic surfaces.
\newblock {\em arXiv preprint arXiv:1202.2304\/} (2012).

\bibitem{IS}
{\sc Ivrissimtzis, I., and Singerman, D.}
\newblock Regular maps and principal congruence subgroups of {H}ecke groups.
\newblock {\em European Journal of Combinatorics 26}, 3 (2005), 437--456.

\bibitem{JS}
{\sc Jones, G.~A., and Singerman, D.}
\newblock Theory of maps on orientable surfaces.
\newblock {\em Proceedings of the London Mathematical Society 3}, 2 (1978),
  273--307.

\bibitem{K}
{\sc Klein, F.}
\newblock \"{U}ber die {T}ransformationen siebenter {O}rdnung der elliptischen
  {F}unktionen.
\newblock {\em Mathematische Annalen 14\/} (1879), 428--471.

\bibitem{S}
{\sc Singerman, D.}
\newblock Universal tessellations.
\newblock {\em Rev. Mat. Univ. Complut. Madrid 1}, 1-3 (1988), 111--123.

\bibitem{SS}
{\sc Singerman, D., and Strudwick, J.}
\newblock Petrie polygons, {F}ibonacci sequences and {F}arey maps.
\newblock {\em Ars Mathematica Contemporanea 10}, 2 (2016), 349--357.

\end{thebibliography}
 \bibliographystyle{acm}

\newpage
\begin{appendix}
\section{Pictures of $\mathcal{M}_3(n)$} For $n=3,\dots, 8$.

\begin{figure}[h!]
\begin{center}
\includegraphics[width=\linewidth]{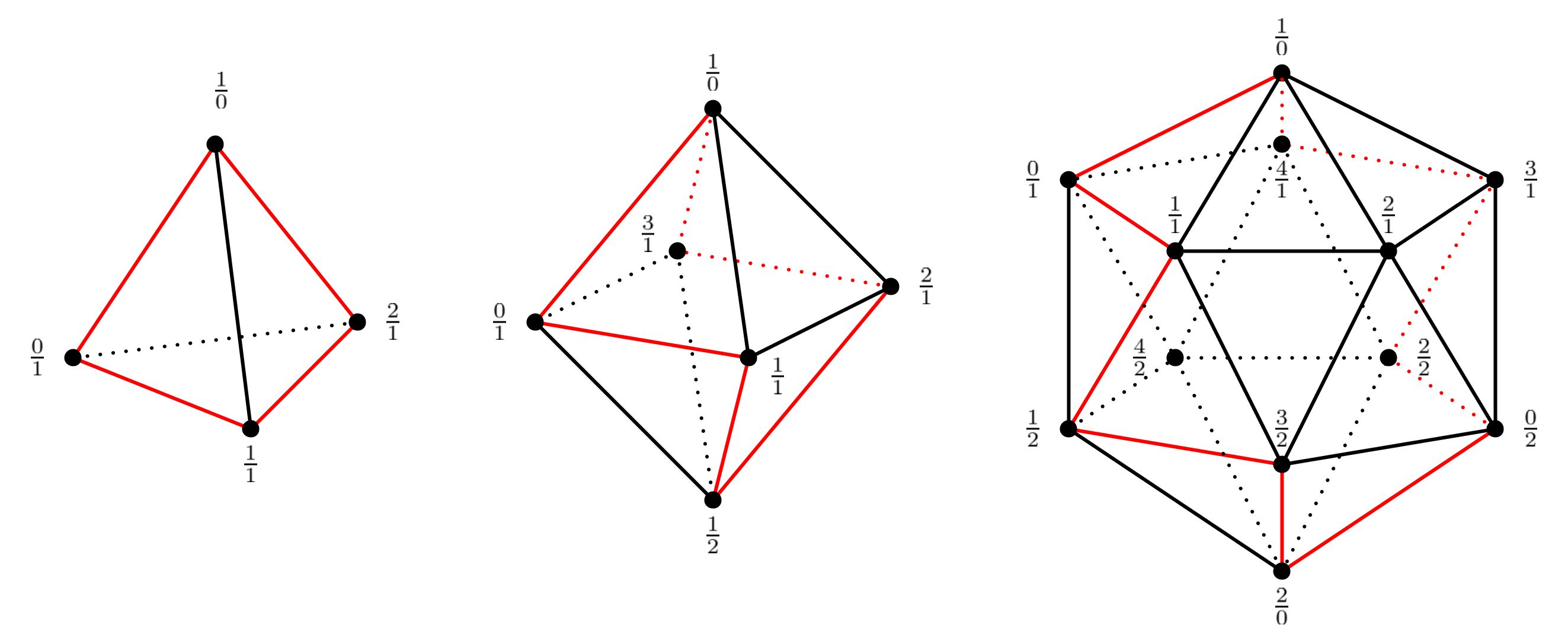}
\end{center}
\caption{Tetrahedron, Octahedron and Icosahedron respectively.}
\label{M345}
\end{figure}
  
\begin{figure}[h!]
\begin{center}
\includegraphics[width=\linewidth]{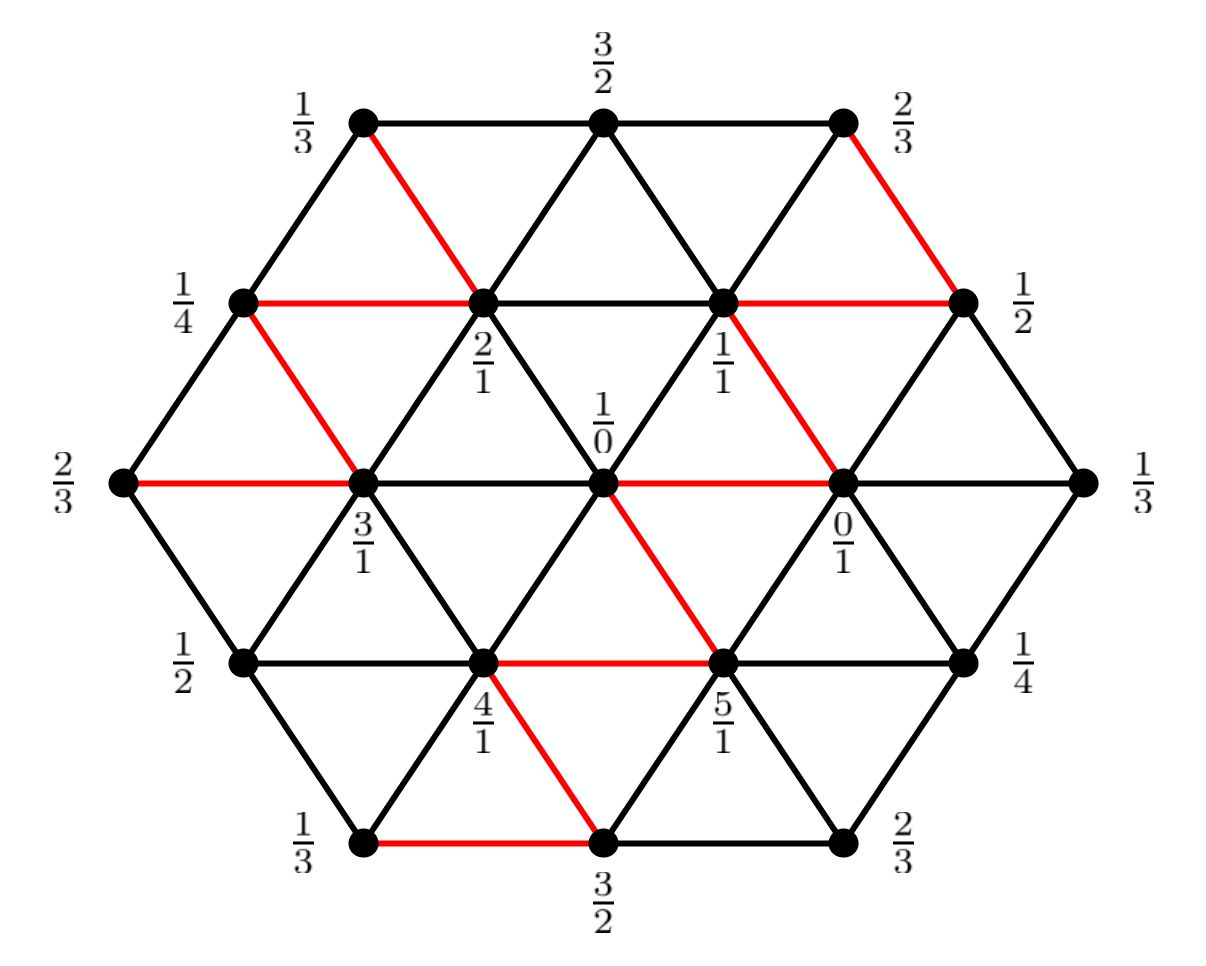}
\end{center}
\caption{$\mathcal{M}_3(6)$}
\end{figure}

\begin{figure}[h!]
\begin{center}
\includegraphics[width=\linewidth]{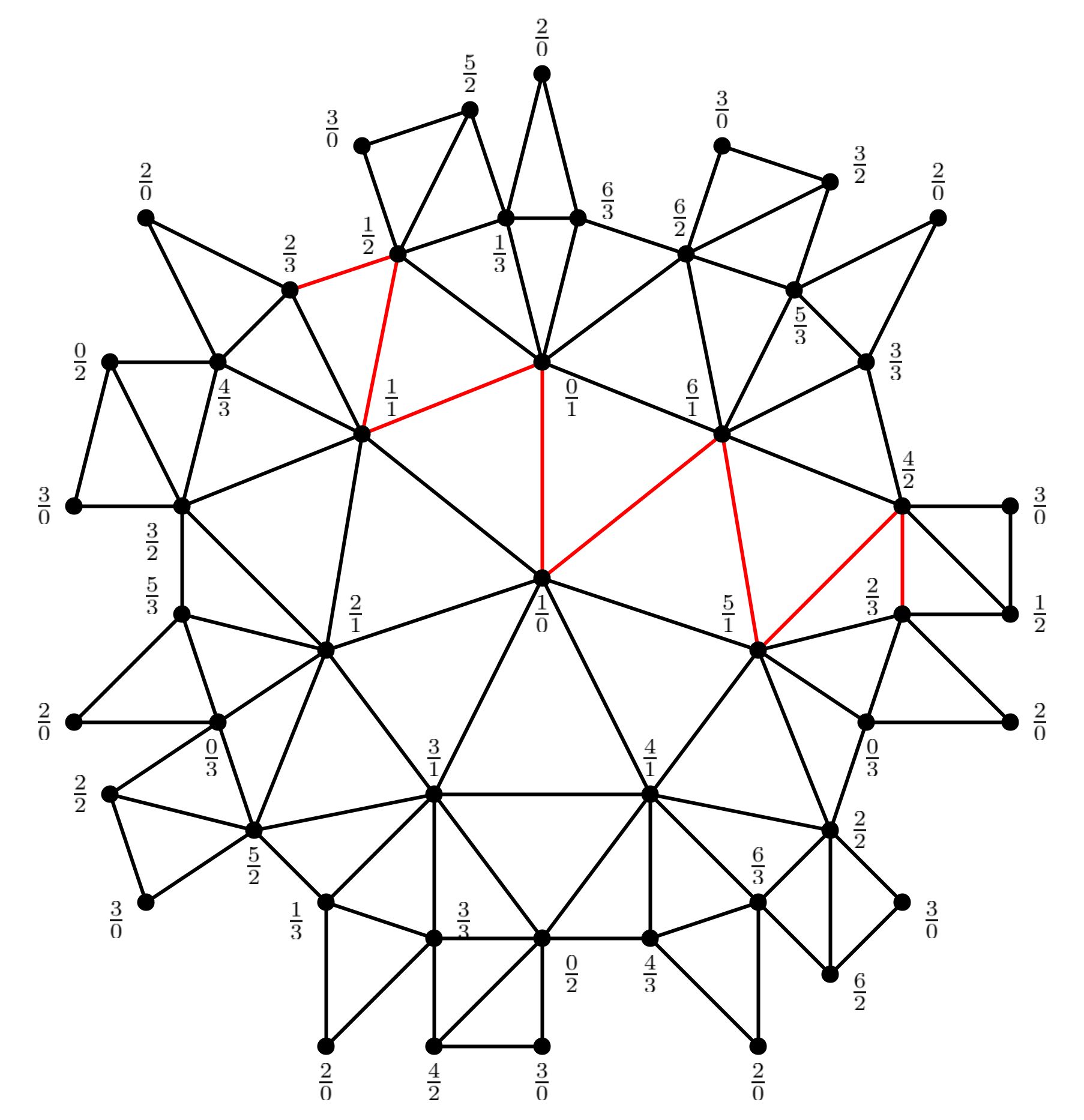}
\end{center}
\caption{$\mathcal{M}_3(7)$ Klein's Riemann surface of genus 3.}
\end{figure}

\begin{figure}[h!]
\begin{center}
\includegraphics[width=\linewidth]{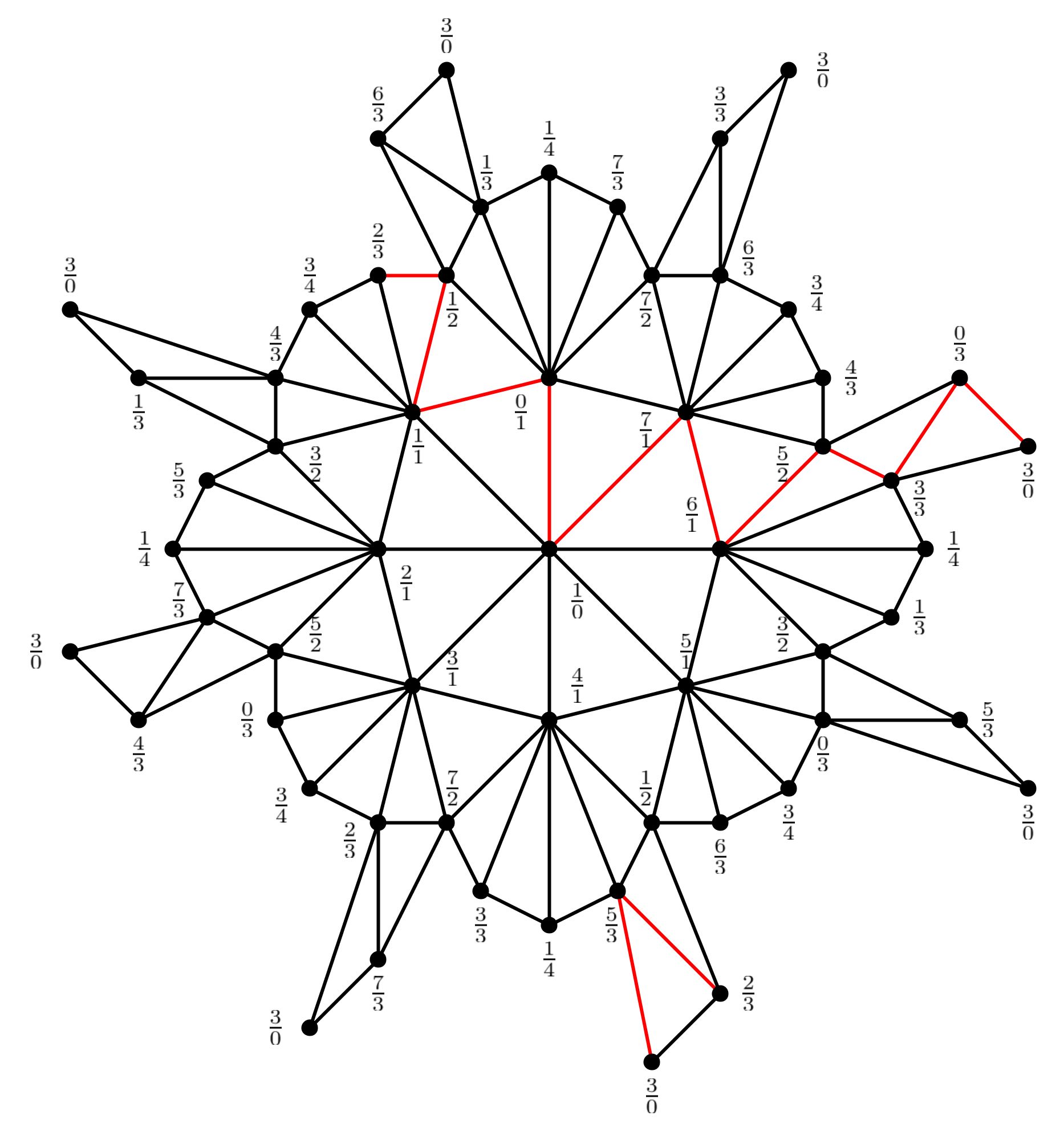}
\end{center}
\caption{$\mathcal{M}_3(8)$} Fricke-Klein surface of genus 5. (Gr\"unbaum polytope.)
\end{figure}

\end{appendix}

\end{document}